\documentclass[12pt,reqno]{amsart}

\usepackage{amssymb, amsmath,}

\renewcommand{\theequation}{\thesection.\arabic{equation}}
\newtheorem{theorem}{Theorem}[section]

\newtheorem{corollary}{Corollary}[section]
\newtheorem{proposition}{Proposition}[section]

\newcommand{\labeq}[1]{\label{eq:#1}}
\newcommand{\refeq}[1]{(\ref{eq:#1})}

\begin{document}

\title{ Generalization of results about the Bohr radius
for power series. }
\author {
Lev Aizenberg}
\thanks{AMS classification number:32A05.}
 \address{
Department of Mathematics \\
           Bar-Ilan University, 52900 Ramat-Gan, Israel}
\email {aizenbrg@math.biu.ac.il}
 \maketitle
 \begin{abstract}
 The Bohr radius for power series of holomorphic functions mapping
  Reinhardt domains ${\mathcal D}\subset {\bf C}^n$ into the convex
domain $G\subset {\bf C}$ is independent of the domain $G.$
 \end{abstract}
 \section{Preliminaries }
\setcounter{equation}{0} Let us recall the theorem of H.Bohr
\cite{bohr:gnus} in 1914.
\begin{theorem} If a power series
\begin{eqnarray}
f(z_1)=\sum\limits_{k=0}^{\infty}c_kz_1^k
\end{eqnarray}
 converges in the
unit disk $U_1$ and its sum has modulus
 less than $1$, then
 \begin{eqnarray}
 \sum\limits_{k=0}^{\infty }\vert c_k z_1^k\vert <1,
 \end{eqnarray}
 if $ \vert z_1\vert <{1\over 3}$. Moreover,
 the constant ${1\over 3}$ cannot be improved.
 \end{theorem}
 For convenience we write the inequality $(1.2)$ in the
 following equivalent form
 \begin{eqnarray}
\sum\limits_{k=1}^{\infty }\vert c_k z_1^k\vert 
<1-\vert c_0\vert . \end{eqnarray}
Later, certain generalizations of this result were obtained.
\newline $\bf{1^0}.$ (\cite{tom:gnus}) If the sum of the
series $(1.1)$ is such that $\vert\Re f(z_1)\vert <1$ in $U_1$ and
$c_0>0$, then for $\vert z_1\vert <{1\over 3}$ the inequality
$(1.3)$ holds.\newline 
$\bf{
2^0}.$(\cite{sid:gnus},\cite{pau:gnus}) If $\Re f(z_1) <1$ in
$U_1$ and $c_0>0$, then for $\vert z_1\vert <{1\over 3}$ the
inequality $(1.3)$ holds.\newline  
$\bf{3^0}.$ (\cite{kres:gnus})
If $\Re \{[\exp(-i arg f(0))]  f(z_1)\} <1$ in $U_1$ (here we assume
that $arg f(0)>0$, if $f(0)=0$), then for $\vert
z_1\vert <{1\over 3}$ the inequality $(1.3)$ is valid.\\

 Formulations  of Bohr's theorem in several complex variables appeared
 very recently. We recall some of them.\\
 Given a complete Reinhardt domain ${\mathcal {D}},$ 
we denote by $R_1({\mathcal D})$ (or by $R_2({\mathcal {D}})$)
 the largest nonnegative number $r$ with he property that if the power series
  \begin{eqnarray}
    f(z)=\sum\limits _{\vert \alpha \vert \geq 0}
    c_{\alpha}z^{\alpha},\;z\in {\mathcal {D}},
    \end{eqnarray}
    where $\alpha =(\alpha _1, \dots , \alpha _n)$,
    $\vert \alpha \vert = \alpha _1+ \dots + \alpha _n $,
    $z^{\alpha}=z_1^{\alpha _1}\dots z_n^{\alpha _n}$ and all $\alpha _i$ are
    nonnegative integers, converges in ${\mathcal {D}}$ and the
 modulus of its sum is less than $1,$ then
 \begin{eqnarray*}
 \sum\limits_{\vert \alpha \vert \geq 1}\vert c_{\alpha } z^{\alpha }
\vert <1-\vert c_0\vert 
 \end{eqnarray*}
 in the homothetic domain  ${\mathcal {D}}_r=r{\mathcal {D}}$. 
Here  $c_0=c_{0,0\dots, 0}.$ Correspondingly, if we consider a
 bounded domain ${\mathcal {D}}$ for $R_2({\mathcal
 D})$ we have
\begin{eqnarray*}
 \sum\limits_{\vert \alpha \vert \geq 1}\sup_{{\mathcal {D}}_r}\vert 
c_{\alpha } z^{\alpha }\vert <1-\vert c_0\vert.
 \end{eqnarray*}
 Let
  \begin{eqnarray*}
 {\mathcal {D}}^n_p=\{z\in {\bf C}^n:\;\vert z_1 \vert^p +\dots +
 \vert z_n\vert^p  < 1\},
 \end{eqnarray*}
 where $0< p \leq \infty $. The domain ${\mathcal {D}}^n_{\infty }$ 
is the poly-disk  $\{z\in {\bf C}^n:\; \vert z_j\vert <1, j=1,\dots, n \}$.
 \begin{theorem} (\cite{boha:gnus}, see also \cite{din:gnus})
 For $n>1 $ one has
 \begin{eqnarray}
 {1\over 3\sqrt n}<R_1({\mathcal {D}}^n_{\infty })<
{2\sqrt{\log n}\over \sqrt n}.
 \end{eqnarray}
 \end{theorem}
\begin{theorem}(\cite{aiz1:gnus})
 For $n>1 $ one has
 \begin{eqnarray}
 {1\over 3\sqrt[3]{e}}<R_1({\mathcal {D}}^n_1)\leq{1\over 3}.
 \end{eqnarray}
 \end{theorem}
 The estimates $(1.5)$ and $(1.6)$ were generalized for $R_1({\mathcal
 D}^n_p)$, for $1\leq p <\infty $ in \cite{boas:gnus} and for $0<p\leq
 1 $ in \cite{aiz3:gnus}. We point out the next new remarkable result
which improves the lower estimate in $(1.5).$
\begin{theorem} (\cite{deff:gnus})
 For $n>1 $ one has
 \begin{eqnarray*}
 C\sqrt{{\log n \over n\log\log  n}}<R_1({\mathcal {D}}^n_{\infty }),
 \end{eqnarray*}
 where the constant $C$ is independent of $n$.
 \end{theorem}
 Both Bohr radii coincide in the case the domain is a polydisk,
 and in the case $n=1$ they do coincide with the classical Bohr
 radius ${1\over 3}$. If the domain ${\mathcal {D}}$ is not a polydisk,
 then naturally $R_2({\mathcal {D}})$ is smaller than $R_1({\mathcal {D}})$.
\begin{theorem} (\cite{aiz1:gnus})
 The inequality
 \begin{eqnarray*}
 1-\sqrt[n]{{2\over 3}}<R_2({\mathcal {D}})
 \end{eqnarray*}
 is true for every complete bounded Reinhardt domain ${\mathcal {D}}$.
 \end{theorem}
 \begin{theorem} (\cite{aiz1:gnus})
 There holds the inequality
 \begin{eqnarray*}
 R_2({\mathcal {D}}^n_1)<{0.44663\over n}.
 \end{eqnarray*}
 \end{theorem}
 The radius $ R_2({\mathcal {D}})$ was a subject of investigation in
 \cite{boas:gnus},\cite{def1:gnus}. Other results about the Bohr
 radius for holomorphic functions can be found in
 \cite{aiz2:gnus},
 \cite{aizad2:gnus},\cite{aizg:gnus},\cite{aizlv:gnus},
\cite{aizvi:gnus},\cite{ben:gnus}, \cite{def:gnus}.

 \section{Generalized Bohr radii}
 \setcounter{equation}{0}
 One of the proofs of Bohr's theorem (Theorem 1.1) is based on the
 Landau inequality \cite{laga:gnus}: if the function $(1.1)$
 satisfies in $U_1$ the inequality $\vert f(z_1)\vert <1 $, then
 $\vert c_k\vert \leq 2(1-\vert c_0\vert)$ holds  for every $k\geq 1 $.
 This inequality can be obtained as a simple consequence of the
 Caratheodory inequality \cite{carath:gnus}: if the function $(1.1)$
 satisfies in $U_1$ the inequality $\Re f(z_1)>0,$ then $\vert
 c_k\vert \leq 2\Re c_0$ is true for every $k\geq 1.$
 Both inequalities are particular cases of a more general assertion.

Let $\tilde G$ be the convex hull of $G.$

 \begin{proposition} (\cite{aiz3:gnus})
 If $f(U_1)\subset G,$  then
 \begin{eqnarray}
 \vert c_k\vert \leq 2dist ( c_0,\partial \tilde G),
 \end{eqnarray}
 for all  $k\geq 1 $.
 \end{proposition}
 Now it is not difficult to prove a generalization of Theorem 1.1.
 Let $G\subset {\bf C}$ be any domain. A point $p\in
 \partial G$ is called a point of convexity if $p\in
 \partial \tilde G.$ A point of convexity $p$ is called regular if there
 exists a disk $U\subset G$ so that $p\in \partial U$.
 \begin{theorem}
 If the function $(1.1)$ is such that $f(U_1)\subset G$,
 with $\tilde G\not={\bf C}$, then for $\vert z_1\vert
 <{1\over 3}$ the inequality
 \begin{eqnarray}
 \sum\limits_{k=1}^{\infty}\vert c_k z_1^k\vert <dist(c_0,
 \partial \tilde G)
 \end{eqnarray}
 is valid. The constant ${1\over 3}$ cannot be improved if
 $\partial G$ contains at least one regular point of
 convexity.
 \end{theorem}
 {\bf Proof:} 1) If $\vert z_1\vert <{1\over 3}$ then 
 $(2.1)$ yields
 \begin{eqnarray*}
 \sum\limits_{k=1}^{\infty}\vert c_k z_1^k\vert <2dist(c_0,
 \partial \tilde G)\sum\limits_{k=1}^{\infty}{1\over 3^k}=dist(c_0,
 \partial \tilde G).
 \end{eqnarray*}
2) We will prove the exactness of the constant  ${1\over 3}$ in
the case the boundary contains at least one regular point of
convexity. In the classical case of Bohr's Theorem 1.1 this
is obtained by considering the family of functions
(\cite{laga:gnus})
\begin{eqnarray}
f(z_1)={\alpha -z_1\over 1-\alpha z_1}, \; 0<\alpha <1.
\end{eqnarray}
Here
\begin{eqnarray*}
 \sum\limits_{k=1}^{\infty}\vert c_k z_1^k\vert =1
\end{eqnarray*}
if and only if $\vert z_1\vert ={1\over 1+2\alpha}$. Furthermore,
taking $\alpha \longrightarrow 1,$ we obtain the desired result.
Note that instead of the family $(2.3)$ one can use the family
$e^{i\phi} f(z_1)$, where $f(z_1)$ is taken from $(2.3)$. In this
case it follows that $c_0=e^{i\phi}\alpha,$ and when $\alpha
\longrightarrow 1$ we get that $c_0$ tends to $\partial U_1$ along
the radius of argument $\phi $. If $G$ is an arbitrary disk $U,$
then, remarking that $(2.2)$ does not change under homotheties and
translations, we deduce the exactness of ${1\over 3}$ in the case
of any disk. Let $\zeta  $ be a regular point of convexity, then
there exists a disk $U\subset G $ such that $\zeta \in(\partial U
)\cap (\partial G)$. Consider the functions $f$ in $(1.1)$ such
that $f(U_1)\subset U$. For suitable $c_0$ (see above) we will
have $dist(c_0, \partial U)=dist(c_0, \partial G)=dist(c_0,
 \partial \tilde G)$. Therefore, in the inequality $(2.2)$
 one cannot take $\vert z_1\vert <r$, where $r>{1\over
 3}$.$ \quad\diamondsuit$\newline
 We remark that Theorem 1.1, the assertion
$\bf{3^0},$ as well as the generalized assertions $\bf{1^0}$ 
and $\bf{2^0}$ are contained in Theorem 2.1.
 For example, in $\bf{1^0}$ no need in assuming $c_0>0,$
 and instead of $(1.3)$ one gets
 \begin{eqnarray*}
 \sum\limits_{k=1}^{\infty}\vert c_k z_1^k\vert <1-\vert \Re
 c_0\vert .
 \end{eqnarray*}
 Similarly in $\bf{2^0}$ no need in assuming $c_0>0,$
 and instead of $(1.3)$ one gets
\begin{eqnarray*}
 \sum\limits_{k=1}^{\infty}\vert c_k z_1^k\vert <1- \Re
 c_0 .
 \end{eqnarray*}
 Let us recall another fact, known earlier:\newline
 $\bf{4^0}$. (\cite{aizad1:gnus}) If $\Re f(z_1)>0$ in $U_1$ and $c_0>0$,
 then for $\vert z_1\vert <{1\over 3}$ the inequality
\begin{eqnarray}
 \sum\limits_{k=1}^{\infty}\vert c_k z_1^k\vert <c_0 
 \end{eqnarray}
holds, and the constant ${1\over 3}$ cannot be improved.
 I thought before that Theorem 1.1. and $\bf{4^0}$ are two different
 facts, having the same Bohr radius. In the light of Theorem 2.1,
 I know now that both results are particular cases of this
 theorem. Now, in the case of $\bf{4^0}$ without the assumption $c_0
 >0,$ we get
 \begin{eqnarray*}
 \sum\limits_{k=1}^{\infty}\vert c_k z_1^k\vert <\Re c_0 
 \end{eqnarray*}
instead of $(2.4).$ 

Theorem 2.1 motivates the following generalization 
of the first and second Bohr radii. Denote by $R_1({\mathcal {D}},G)$ 
(or by $R_2({\mathcal {D}},G)$), where $G\subset
 {\bf C}$, $\tilde G\not= {\bf C},$ and ${\mathcal {D}}$ is a complete
 Reinhardt domain (bounded complete Reinhardt domain) in ${\bf
 C}^n$ the largest $r\geq 0$ such that if the function $(1.4)$
 is holomorphic in ${\mathcal {D}}$ and $f({\mathcal {D}})\subset G$ then
 \begin{eqnarray*}
 \sum\limits_{\vert\alpha \vert \geq 1}^{\infty}\vert c_{\alpha} 
z^{\alpha}\vert <dist(c_0 ,\partial \tilde G)
 \end{eqnarray*}
 in a homothety ${\mathcal {D}}_r$  (or correspondingly
  \begin{eqnarray*}
 \sum\limits_{\vert\alpha \vert \geq 1}^{\infty}\sup_{{\mathcal 
{D}}_r}\vert c_{\alpha} z^{\alpha}\vert <dist(c_0 ,\partial \tilde G)).
 \end{eqnarray*}
 Theorem 2.1 and the result from \cite{aize:gnus} about the Rogosinski
 radius allow one to hope that the two Bohr radii $R_1({\mathcal {D}},G)$ and 
$R_2({\mathcal D},G)$ are independent of the convex domain $G$. The main
 result of the present paper is the proof of the validity of this more
 general assertion.

 \section{ The main result}
 \setcounter{equation}{0}
 Let $M$ be a complex manifold, ${\mathcal H}(M)$ be the
 space of holomorphic on $M$ functions equipped with the
 natural topology of uniform
 convergence over compact subsets of $M$.\\
 Let $\| \cdot\|_r$, $r\in (0,1)$, be a one-parameter family of
 semi-norms in ${\mathcal H}(M)$ that are continuous with
 respect to the topology of ${\mathcal H}(M)$. In what follows we 
always assume that
 \begin{eqnarray*}
 a)\;\;\| \cdot\|_{r_1}&\leq &\| \cdot\|_{r_2}\;\;if \;r_1\leq r_2.\\
 b)\;\;\| f\cdot g\|_{r}&\leq &\| f\|_{r}\cdot\| g\|_{r}\;\;\forall
 r\in (0,1).
 \end{eqnarray*}
 There exists a point $z_0\in M$ such that
  \begin{eqnarray*}
 c)\;\;\| f\|_{r}&\longrightarrow & \vert f(z_0)\vert \;\rm{as}\; 
r\longrightarrow 0,\; \forall f\in {\mathcal {H}}(M).\\
 d)\;\;\| f\|_{r}&= &\vert f(z_0)\vert+\| f-f(z_0)\|_{r},\;\forall
 f\in {\mathcal {H}}(M).
 \end{eqnarray*}
 Denote by $B(\| \cdot\|_{r}, G)$ the largest $r\geq 0$ such
 that for $f\in {\mathcal {H}}(M)$ and $ f(M)\subset  G$ one has
 \begin{eqnarray}
\| f-f(z_0)\|_{r}<dist (f(z_0), \partial \tilde G),
\end{eqnarray}
where $\tilde G$ is the convex hull of the domain $G\subset {\bf
C}$.
\begin{proposition}
If $U$ is any disk and $\Pi $ is any half-plane, then
\begin{eqnarray}
B(\| \cdot\|_{r}, \Pi)=B(\| \cdot\|_r, U).
\end{eqnarray}
\end{proposition}
{\bf {Proof:}} Let $\Pi _1=\{z_1:\; \Re z_1>0 \}$, then (\cite
{aizad1:gnus}, Theorem 7) 
\begin{eqnarray*}
B(\| \cdot\|_r, U_1)=B^{\prime}(\| \cdot\|_{r}, \Pi _1),
\end{eqnarray*}
where $B^{\prime}$ is defined in the same way as $B$ but with the
additional assumption $f(z_0)>0$. This assumption can be removed
as follows. If $\Re f(z_0)>0$ in $M$ then $\Re f_1(z_0)>0$, where
$f_1(z)=f(z)-\Im f(z_0)$. But $f_1(z_0)>0$, hence
\begin{eqnarray*}
B(\| \cdot\|_r, U_1)=B(\| \cdot\|_{r}, \Pi _1).
\end{eqnarray*}
We remark that $(3.1)$ does not change under homotheties,
translations and rotations of the domain $G$. Therefore
$(3.2)$ holds. $\diamondsuit $\\
\begin{theorem}
If $ \tilde G\not={\bf C}$, then $B(\| \cdot\|_{r}, G)$  is not
smaller than $(3.2)$. If $\partial G$ contains at least one regular point
of convexity, then  $B(\| \cdot\|_{r}, G)$ is equal to
$(3.2)$.
\end{theorem}
{\bf {Proof:}} Let $ \tilde G\not={\bf C}$ and $ f(M)\subset G$.
Fix any $f(z_0)\in G$. On the boundary $\partial \tilde G$ there
exists a point $\zeta $ so that $dist(f(z_0),\partial \tilde G)=
dist(f(z_0), \zeta)$. Through the point $\zeta $ the line
of support of $\tilde G$ passes which defines the half-plane $\Pi _0\supseteq
G$. Then $$dist(f(z_0),\partial \tilde G)=dist(f(z_0),\partial
\Pi_0).$$ Therefore $B(\| \cdot\|_{r}, G)\geq B(\| \cdot\|_{r}, \Pi
_0)$, since $\{f:f\in {\mathcal {H}}(M),\; f(M)\subset G\}\subset
\{f:f\in {\mathcal {H}}(M),\; f(M)\subset \Pi _0\}$.\\
Assume now that there is a regular point of convexity in $\partial
G$. Then the proof repeats the proof of part 2) of the Theorem
2.1. Note that there we did not use the concrete form of the
family $(2.3)$, but rather the fact that $c_0$ can lie
on any radius emanating from the center of the disk $U$ to its
boundary. So, let us assume that $U\subset G$, $\zeta \in
(\partial U)\cap (\partial G)\cap (\partial\tilde G)$. Then
consider $f(z_0)$ lying on the radius from the center of the disk
$U$ to the point $\zeta $. Now $dist(f(z_0),
\zeta)=dist(f(z_0),\partial U)=dist(f(z_0),\partial \tilde G)$,
hence  $B(\| \cdot\|_{r}, G)\leq B(\| \cdot\|_{r}, U)$, since
$\{f:f\in {\mathcal {H}}(M),\; f(M)\subset U\}\subset \{f:f\in
{\mathcal {H}}(M),\; f(M)\subset G\}$. $\diamondsuit $.
\begin{corollary}
If the domain $G$ is convex and $G\not={\bf C}$, then $B(\|
\cdot\|_{r}, G)$ is independent of the choice of the domain $G$.
\end{corollary}
{\bf {Proof:}} There exists disk $U\subset G$ such that $\partial
U\cap \partial  G\not=\emptyset $. Therefore there exist regular
points of convexity on $\partial  G$.
$\diamondsuit $\\
\begin{corollary}
The first Bohr radius $R_1({\mathcal {D}}, G)$ and the second Bohr
radius $R_2({\mathcal {D}}, G)$ are independent of the choice of
the convex domain $G$, $G\not={\bf C}$.
\end{corollary}
In particular, the assertions of Theorems 1.2, 1.3 and 1.4 are valid
for $R_1({\mathcal {D}}, G)$ while those of Theorem 1.5 and 1.6
are valid for $R_2({\mathcal {D}}, G)$ for every convex domain
$G\not={\bf C}$.\\

{\bf Some concluding remarks.} If the family of semi-norms $\|\cdot\|_r
$ does not satisfy some of the conditions $a)-d)$, then the assertion
of Theorem 3.1 is not valid anymore. Examples can be found in
\cite{aizad1:gnus}. If $\tilde G={\bf C}$, then the right-hand 
side of $(3.1)$ is equal to $\infty $, therefore in this case
$$B(\| \cdot\|_{r}, G)=1.$$ One can also consider different
realizations of $B(\| \cdot\|_{r}, G)$ than the first and second
Bohr radii $R_1({\mathcal {D}}, G)$ and $R_2({\mathcal {D}},
G)$.\\
We conclude the present article with formulating an {\bf{open
problem}}: if $\tilde G\not={\bf C},$ is it always true that
$B(\| \cdot\|_{r}, G)$ is equal to $(3.2)$? The same question makes
sense for the first and second Bohr radii $R_1({\mathcal {D}},
G)$ and $R_2({\mathcal {D}}, G)$.

\section{Acknowledgements}

The author is deeply grateful to E. Liflyand and A. Vidras for their
help in preparing the paper and improving the presentation.

\end{document}